\documentclass{amsart}
\usepackage[dvips,final]{epsfig}
\usepackage{amsmath}
\usepackage{amscd}
\usepackage{amsfonts}
\usepackage{amssymb}
\usepackage{latexsym}
\newtheorem{theorem}{Theorem}%[section]
\newtheorem{lemma}{Lemma}
\newtheorem{corollary}{Corollary}

\theoremstyle{definition}

\theoremstyle{remark}

\numberwithin{equation}{section}

\begin{document}

\title{Multicusps}


%    author one information
 \author[Y.~Mizota and T.~Nishimura]{Yusuke Mizota and Takashi Nishimura}
%\author{Yusuke Mizota}
\address{Graduate School of Mathematical Sciences, 
Kyushu University, 
744, Motooka, Nishi-ku, Fukuoka 819-0395, JAPAN. 
}
%\curraddr{}
\email{ma210046@math.kyushu-u.ac.jp}
%\thanks{}

%    author two information
%\author{Takashi Nishimura$^*$}
%\thanks{$^*$Corresponding author. }
\address{Research Group of Mathematical Sciences,  
Research Institute of Environment and Information Sciences, 
%Faculty of Education and Human Sciences, \\
Yokohama National University, 
%	79-2 Tokiwadai, Hodogaya-Ku, \\
Yokohama 240-8501, JAPAN.  
}
%\curraddr{}
\email{nishimura-takashi-yx@ynu.jp}
%\subjclass[2000]{Primary }
%    The 2010 edition of the Mathematics Subject Classification is
%    now available.  If you are citing a classification from the
%    new scheme, use the following input coding instead.
\subjclass[2010]{Primary 57R45; Secondary 58K20, 58K45. }
\keywords{Multicusp, reduced Kodaira-Spencer-Mather map, 
higher version of the reduced Kodaira-Spencer-Mather map, liftable vector field}

\date{}

\begin{abstract}
For a given multicusp $f=c_{(\theta_0, \ldots, \theta_i)}$ $(1\le i)$,  
we present a direct sum decomposition theorem 
of the source space of ${}_i\overline{\omega}f$, 
where ${}_i\overline{\omega}f$ is 
a higher version of  the reduced Kodaira-Spencer-Mather map $\overline{\omega}f$.   
As a corollary of our direct sum decomposition theorem, we show that 
for any $i\in \mathbb{N}$ and any  $f=c_{(\theta_0, \ldots, \theta_i)}$,  
${}_i\overline{\omega}f$ is bijective.    
The corollary is an affirmative answer to the question raised by M. A. S. Ruas 
during the 11th International Workshop on Real and Complex Singularities 
at the University of S${\tilde {\rm a}}$o Paulo in S${\tilde {\rm a}}$o Carlos (2010).    
\end{abstract}

\maketitle
\section{Introduction} % use lowercase except for proper names
\label{intro}
%\noindent 
Throughout this paper, let $\mathbb{K}$ be $\mathbb{R}$ or $\mathbb{C}$, and let $S$ be a finite subset of $\mathbb{K}^n$.   
% and $f: (\mathbb{K}^n,S) \to (\mathbb{K}^p,0)$ $(n\le p)$ be an 
%For a finite subset $S$ of $\mathbb{K}^n$ and $0$ of  $\mathbb{K}^n$,  
Let $C_S$ (resp., $C_0$) be the $\mathbb{K}$-algebra of $C^\infty$ or holomorphic function germs $(\mathbb{K}^n,S)\to \mathbb{K}$ 
(resp., $(\mathbb{K}^p, 0)\to \mathbb{K}$), and let $m_S$ (resp., $m_0$) be the subset of 
$C_S$ (resp., $C_0$) consisting of $C^\infty$ or holomorphic function germs $(\mathbb{K}^n,S)\to (\mathbb{K},0)$ 
(resp., $(\mathbb{K}^p, 0)\to (\mathbb{K},0)$).        
%The sets $C_S$ and $C_0$ have natural $\mathbb{K}$-algebra structures induced 
%by the $\mathbb{K}$-algebra structure of $\mathbb{K}$.   
A map germ $f: (\mathbb{K}^n,S)\to (\mathbb{K}^p,0)$ is called a {\it multigerm}.   For a $C^\infty$ or holomorphic multigerm $f: (\mathbb{K}^n,S)\to (\mathbb{K}^p,0)$, we let $f^*: C_0\to C_S$ be the 
$\mathbb{K}$-algebra homomorphism defined by $f^*(u)=u\circ f$.   
%Put 
%$Q(f)=C_S/f^*m_0C_S$.      
%An smooth multigerm  $f: (\mathbb{K}^n,S) \to (\mathbb{K}^p,0)$ is said to be {\it finite} 
%if $\dim_{\mathbb{K}}Q(f)$ is finite.      It is known that if a mulltigerm  $f: (\mathbb{K}^n,S) \to (\mathbb{K}^p,0)$ is finite, then 
%$n\le p$ and in the case that $n\le p$ finite multigerms  $f: (\mathbb{K}^n,S) \to (\mathbb{K}^p,0)$ consist of 
%a generic set (for instance, see \cite{wall}).   
\par 
For a $C^\infty$ or holomorphic map germ 
%$C^\infty$ map-germ 
$f: (\mathbb{K}^n, S)\to \mathbb{K}^p$,  
%such that $f(S)\subset T$, 
%where $T$ 
%is a finite subset 
%of $\mathbb{K}^p$, 
%be a germ of $C^\infty$ mapping such that 
%$f(S)\subset T$. 
%For the $f$ 
let $\theta_S(f)$ be the $C_S$-module consisting of germs of 
$C^\infty$ or holomorphic vector fields along $f$.   
%If $f$ is a multigerm $f: (\mathbb{K}^n, S)\to (\mathbb{K}^p,0)$
The $C_S$-module $\theta_S(f)$ is naturally identified with 
$\underbrace{C_S\oplus \cdots \oplus  
C_S}_{p \mbox{ tuples}}$.    
We define $\theta_S(n)$ (resp., $\theta_{0}(p)$) 
as the $C_S$-module (resp., $C_0$-module) 
$\theta_S(id._{({\mathbb{K}^n},S)})$ (resp., $\theta_{\{0\}}(id._{(\mathbb{K}^p,0)})$, 
where $id._{({\mathbb{K}^n},S)}$ 
(resp., $id._{(\mathbb{K}^p,0)}$) is the germ of the identity map of  
%map-germ of 
$(\mathbb{K}^n,S)$ (resp., $(\mathbb{K}^p,0)$).   
%and 
%$id._{(\mathbb{K}^p,0)}$ is the identity map-germ of $(\mathbb{K}^p,0)$.    
%, 
%(If $k=\infty$, then $k-1$ is again $\infty$.   
%That is to say, $m_q^\infty\theta(q)$ is the set of all flat map-germs $(q=n, p%))$.   
%Since $m_q^\aleph=\{0\}$, $m_q^\aleph\theta(q)=\{0\}$ $(q=n, p)$.)     
\par 
For a given $C^\infty$ or holomorphic multigerm 
$f: (\mathbb{K}^n,S)\to (\mathbb{K}^p,0)$, 
%the $\mathcal{G}$-equivalence class of $f$ is denoted by $\mathcal{G}(f)$, 
%where $\mathcal{G}$ is one of $\mathcal{L, A, C}$ and $\mathcal{K}$.   
%For the $f$, 
%following Mather (\cite{mather3}), we 
define $tf: \theta_S(n)\to \theta_S(f)$ and  $\omega f : \theta_{0}(p) \to \theta_S(f)$ as 
%follows: 
%\begin{eqnarray*}
$tf(\eta)=df\circ \eta$ and $\omega f(\xi)= \xi\circ f$,  respectively,     
%\end{eqnarray*}
where $df$ is the differential of $f$ (for details on $tf$ and $\omega f$, see \cite{mather3}).   
For the $f$, put $T\mathcal{R}_e(f)  =  tf(\theta_S(n))$, $T\mathcal{L}_e(f)  =  \omega f(\theta_{0}(p))$, 
and $T\mathcal{A}_e(f)  =  T\mathcal{R}_e(f)+T\mathcal{L}_e(f)$ 
(for details on  $T\mathcal{R}_e(f), T\mathcal{L}_e(f)$, 
and $T\mathcal{A}_e(f)$, see \cite{wall}).    
The set of $C^\infty$ or holomorphic function germs $(\mathbb{K}^n,S)\to \mathbb{K}$ 
(resp., $(\mathbb{K}^p,0)\to \mathbb{K}$), such that 
the terms of their Taylor series up to $(i-1)$ are zero, is denoted by $m_S^i$ 
(resp., $m_0^i$) for any non-negative integer $i$.     
Therefore, $m_S^0$ (resp., $m_0^0$) is exactly the same as $C_S$ (resp., $C_0)$.     
A $C^\infty$ or holomorphic multigerm  
$f: (\mathbb{K}^n,S)\to (\mathbb{K}^p,0)$ is said to be {\it finitely determined}  if there exists a positive integer 
$k$ such that the inclusion $m_S^k\theta_S(f)\subset T\mathcal{A}_e(f)$ holds.   
%It is known that any finitely determined multigerm $f: (\mathbb{K}^n,S)\to (\mathbb{K}^p,0)$ $(n\le p)$ has the property that 
%$\dim_{\mathbb{K}}Q(f)$ is finite (for instance, see \cite{wall}).   
%\medskip    
\par 
For a given $C^\infty$ or holomorphic multigerm $f: (\mathbb{K}^n, S)\to (\mathbb{K}^p,0)$, 
a vector field $\xi\in \theta_0(p)$ is said to be 
{\it liftable over $f$}  if $\xi\circ f$ belongs to $T\mathcal{L}_e(f)\cap T\mathcal{R}_e(f)$.    The set of vector fields liftable over $f$ naturally has a $C_0$-module structure.      
In \cite{nishimuraliftable}, in order to express the minimal number of generators 
for the module of vector fields liftable over $f$, 
the {second} author introduced 
the following homomorphism ${}_i\overline{\omega}f$ for a  given $C^\infty$ or holomorphic multigerm $f$ and 
a non-negative integer $i$.   
\begin{eqnarray*}
{}_i\overline{\omega}f : \frac{m_0^i\theta_0(p)}{m_0^{i+1}\theta_0(p)}
& \to &  
\frac{f^*m_0^i\theta_S(f)}{T\mathcal{R}_e(f)\cap f^*m_0^i\theta_S(f)+f^*m_0^{i+1}\theta_S(f)},  \\
{}_i\overline{\omega}f([\xi])
& = & 
[\omega f(\xi)].   
\end{eqnarray*}
Note that ${}_0\overline{\omega}f$ is identical to the map $\overline{\omega}f$ defined in \cite{mather4}, 
i.e.,    
\begin{eqnarray*}
\overline{\omega}f : \frac{\theta_0(p)}{m_0\theta_0(p)}
& \to &  
\frac{\theta_S(f)}{T\mathcal{R}_e(f)+f^*m_0\theta_S(f)},  \\
\overline{\omega}f([\xi])
& = & 
[\omega f(\xi)].   
\end{eqnarray*}
Suppose that $\mathbb{K}=\mathbb{C}, n\ge p$ and $S=\{\mbox{one point }x\}$.    
Then, the map $\overline{\omega}f$ is called the {\it reduced Kodaira-Spencer map}, 
and it is denoted by 
$\rho_f(x)$ in \cite{looijenga}.     
In this paper, we are mainly interested in the case $n < p$ and 
$S$ is not a set of only one point; hence, we call $\overline{\omega}f$ (resp., ${}_i\overline{\omega}f$) 
the {\it reduced Kodaira-Spencer-Mather map} 
(resp., a {\it higher version of the reduced Kodaira-Spencer-Mather map}).   
When $n \le p$, the module of vector fields liftable over $f$ can be investigated  
by using higher versions of the reduced Kodaira-Spencer-Mather map as follows.          
%The homomorphism ${}_i\overline{\omega}f$ is it is clearly seen that ${}_i\overline{\omega}f$ is well-defined.       
%Note that ${}_0\overline{\omega}f=\overline{\omega}f$.  
\begin{theorem}[\cite{nishimuraliftable}]
Let $f: (\mathbb{K}^n, S)\to (\mathbb{K}^p,0)$ be a finitely determined multigerm of corank at most one.     
Suppose that there exists a non-negative integer 
%such that $\delta(f)< \infty$ 
$i$ such that ${}_i\overline{\omega}f$ is bijective.   
Then, the minimal number of generators for the module of vector fields liftable over $f$ is 
exactly $\dim_{\mathbb{K}}\mbox{\rm ker}({}_{i+1}\overline{\omega}f)$.  
%\dim_{\mathbb{K}}\mbox{\rm ker}({}_{i+1}\overline{\omega}f)=
\end{theorem}
Here, 
{\it corank at most one} for a multigerm $f: (\mathbb{K}^n, S)\to (\mathbb{K}^p,0)$ 
implies that $n\le p$ and $\max\{n-\mbox{rank} Jf(s_j)\; |\; 1\le j\le |S|\}\le 1 $ holds, where 
$Jf(s_j)$ is the Jacobian matrix of $f$ at $s_j\in S$.     
\par 
\medskip 
In \cite{nishimuraliftable}, we may find several examples satisfying the assumption of Theorem 1.     
%However, 
Unfortunately,  there are no examples satisfying the condition that ${}_i\overline{\omega}f$ is bijective 
for some $i$ such that $2\le i$.      
In this paper, for any $i\in \mathbb{N}$, we give a concrete multigerm $f$ such that 
${}_i\overline{\omega}f$ is bijective.     
%The existence of such examples was questioned by 
M.~A.~S.~Ruas asked for such examples 
during the 11th International Workshop on Real and Complex Singularities 
at the University of S${\tilde {\rm a}}$o Paulo in S${\tilde {\rm a}}$o Carlos (2010).        
Thus, this paper answers her question affirmatively.        
\par 
\medskip 
Let $c: \mathbb{K}\to \mathbb{K}^2$ be the map defined by 
$c(x)=(x^2, x^3)$, and for any real number $\theta$, let 
$R_{\theta}: \mathbb{K}^2\to \mathbb{K}^2$ be the linear map that gives 
the rotation of $\mathbb{K}^2$ about the origin with respect to the angle $\theta$.     
$$
R_{\theta}
\left(
\begin{array}{c}
X \\ 
Y
\end{array}
\right)
= 
\left(
\begin{array}{rr}
\cos\theta & -\sin\theta \\
\sin\theta & \cos\theta
\end{array}
\right)
\left(
\begin{array}{c}
X \\ 
Y
\end{array}
\right)
.
$$   
Let $\theta_0, \ldots, \theta_i$ be real numbers satisfying $0\le \theta_j<2\pi$ $(0\le j\le i)$ and 
$0\ne |\theta_j-\theta_k|\ne \pi$ $(j\ne k)$.     
Put $S=\{s_0, \ldots, s_i\}$ $(s_j\ne s_k \mbox{ if }j\ne k)$ and define 
$c_{\theta_j}: (\mathbb{K}, s_j)\to (\mathbb{K}^2,0)$ as $c_{\theta_j}(x)=R_{\theta_j}\circ c(x_j)$, 
where $x_j=x-s_j$.       
%Then, the set $\{c_{\theta_0}, \ldots, c_{\theta_i}\}$ may be regarded as a multigerm 
%$(\mathbb{K},S)\to (\mathbb{K}^2,0)$ defined by $c_{\theta_j}(x_j)$ $(0\le j\le i)$ where $x_j=x-s_j$.      
A multigerm $\{c_{\theta_0}, \ldots, c_{\theta_i}\}: (\mathbb{K},S)\to (\mathbb{K}^2,0)$ 
is called a {\it multicusp}, and  
%A multigerm  $c_{(\theta_0, \ldots, \theta_{i})} : (\mathbb{K},S)\to (\mathbb{K}^2,0)$  
it is denoted by 
$c_{(\theta_0, \ldots, \theta_i)}$.     
In the case $i=0$ (resp., $i=1$), it is called a {\it cusp} (resp., {\it double cusp}).         
It is known that double cusps are open $\mathcal{A}$-simple multigerms $(\mathbb{K},S)\to (\mathbb{K}^2,0)$ 
(see \cite{kolgushkinsadykov}), 
and to the best of authors' knowledge, there is no literature on multicusps for the case 
$2\le i$.       
%\par  
%\medskip 
\par 
For a given multicusp $c_{(\theta_0, \ldots, \theta_{i})}$ such that $i\in \mathbb{N}$,  
put 
$$
c_{\hat{\theta_j}}=\{c_{\theta_0}, \ldots, {\hat{c_{\theta_j}}}, \ldots, c_{\theta_i}\}
: (\mathbb{K},S-\{s_j\})\to (\mathbb{K}^2,0), 
$$ 
where 
${\hat{c_{\theta_j}}}$ denotes the removal of the cusp ${c_{\theta_j}}$.     Note that in order to define $c_{\hat{\theta_j}}$, 
$i$ must be positive.     The main result of this paper can be stated 
in terms of ${}_i\overline{\omega}c_{\hat{\theta_j}}$ as follows.   
%%%%%%%%%%%%%%%%%%%%%%%%%%%%%%%%%%%%%%%%%%%%%%%%%%%%%%%%%%    
\begin{theorem}
For any $i\in \mathbb{N}$ and any multicusp $c_{(\theta_0, \ldots, \theta_{i})}$,  
the following equality holds.  
$$
\frac{m_0^i\theta_0(2)}{m_0^{i+1}\theta_0(2)}=
\bigoplus_{j=0}^i\mbox{\rm ker}({}_i\overline{\omega}c_{\hat{\theta_j}}).   
$$
\end{theorem}
%In the case that $i=0$ Theorem 2 does not hold since ${}_0\overline{\omega}c_{\theta}$ is not the zero map for %any cusp $c_{\theta}$.    
\begin{corollary}
For any $i\in \mathbb{N}$ and any multicusp $c_{(\theta_0, \ldots, \theta_{i})}$, ${}_i\overline{\omega}c_{(\theta_0,\ldots, \theta_i)}$ 
is bijective.   
\end{corollary}
Note that the reduced Kodaira-Spencer-Mather map $\overline{\omega}c_{\theta}={}_0\overline{\omega}c_{\theta}$ 
is not surjective for any cusp $c_{\theta}$ \cite{mather4} 
because $\theta_S(f)\ne T\mathcal{A}_e(f)$ for $f=c_\theta$ and $S=\{0\}$.   
%$c_{\theta}$ is not stable.     
Thus, Corollary 1 does not hold in the case $i=0$.       
Since the image of $c_{(\theta_0, \ldots, \theta_i)}$ is a plane curve, the minimal number of 
generators for the module of vector fields liftable over 
$c_{(\theta_0, \ldots, \theta_i)}$ is always $2$ in the  complex case \cite{saito}.   
This implies that for any $j$ $(0\le j\le i)$, 
$$
\dim_{\mathbb{C}}\left(\mbox{\rm ker}\left({}_i\overline{\omega}c_{\hat{\theta_j}}\right)\right)=2.
$$
Nevertheless, since we want to study both real and complex cases in this paper, 
%without using \cite{saito} 
we directly show that for any $j$ $(0\le j\le i)$, 
$$
\dim_{\mathbb{K}}\left(\mbox{\rm ker}\left({}_i\overline{\omega}c_{\hat{\theta_j}}\right)\right)=2
$$
by using linear algebra (see Lemma 1 in \S 2).     
\par 
For any $i\in \mathbb{N}$, let $j$ be an integer such that $0\le j\le i$.    
Then, since the kernel of 
${}_i\overline{\omega}c_{\hat{\theta_j}}$ is a $2$-dimensional vector space, 
%by Theorem 2, 
there exists a non-zero vector field $\xi_j\in m_0^i\theta_0(2)$ such that 
$\omega c_{\hat{\theta}_j}(\xi_j)\in T\mathcal{R}_e(c_{\hat{\theta}_j})+c_{\hat{\theta}_j}^*m_0^{i+1}\theta_S(c_{\hat{\theta}_j})$.    
The map ${}_{i+1}\overline{\omega}c_{\hat{\theta}_j}$ is surjective 
because the map ${}_{i}\overline{\omega}c_{\hat{\theta}_j}$ is surjective by Corollary 1.     
%and hence 
%The map ${}_{i}\overline{\omega}c_{\hat{\theta}_j}$ is surjective.       
Thus, by the preparation theorem, there 
must exist a vector field $\widetilde{\xi}_j\in m_0^{i+1}\theta_0(2)$ such that 
$\omega c_{\hat{\theta}_j}(\xi_j+\widetilde{\xi}_j)\in T\mathcal{R}_e(c_{\hat{\theta}_j})$.       
On the other hand, $\omega c_{\theta_j}(\xi_j+\widetilde{\xi}_j)\not \in T\mathcal{R}_e(c_{\theta_j})$ by Theorem~2.     
Therefore, by integrating $\xi_j+\widetilde{\xi}_j+\frac{\partial}{\partial t}$, we obtain a diffeomorphism germ    
%By combining Theorem 2 and the proof of Theorem 1, we see that for any $j$ $(0\le j\le i)$ 
%there exists a germ of diffeomorphism 
$H_j: (\mathbb{K}^2,0)\to (\mathbb{K}^2,0)$ 
satisfying the following 4 conditions:
\begin{enumerate}
\item The ${(i-1)}$-jet of $H_j$ at the origin is equal to the $(i-1)$-jet of 
$id._{(\mathbb{K}^2,0)}$ at the origin.      
\item The $i$-jet of $H_j$ at the origin is slightly different from the $i$-jet of $id._{(\mathbb{K}^2,0)}$ at the origin.          
\item The image of $H_j\circ c_{\hat{\theta_j}}$ is equal to 
the image of  $c_{\hat{\theta_j}}$.            
\item The image of $H_j\circ c_{\theta_j}$ is not equal to the image of  $c_{\theta_j}$.     
\end{enumerate}           
\par 
\noindent 
Similarly, 
%Furthermore, by Corollary 1 and the proof of Theorem 1, 
we obtain a diffeomorphism germ 
$\widetilde{H}: (\mathbb{K}^2,0)\to (\mathbb{K}^2,0)$  
satisfying the following 4 conditions:
\begin{enumerate}
\item The ${i}$-jet of $\widetilde{H}$ at the origin is equal to the ${i}$-jet of $id._{(\mathbb{K}^2,0)}$ at the origin.     
\item The ${(i+1)}$-jet of $\widetilde{H}$ at the origin is slightly different from the ${(i+1)}$-jet of $id._{(\mathbb{K}^2,0)}$ 
at the origin.           
\item The image of $\widetilde{H}\circ c_{(\theta_0,\ldots, \theta_i)}$ is equal to 
the image of  $c_{(\theta_0,\ldots, \theta_i)}$.           
\item There exists a sufficiently small $x\in \mathbb{K}$ such that 
$\widetilde{H}\circ c_{(\theta_0, \ldots, \theta_i)}(x)\ne c_{(\theta_0, \ldots, \theta_i)}(x)$.    
\end{enumerate}
\par 
\medskip
\noindent 
\underline{\it Proof of Corollary 1.}\qquad     
Since $1\le i$, by Theorem 2, we have  
$$
\mbox{\rm ker}({}_i\overline{\omega}c_{(\theta_0, \ldots, \theta_i)})
\subset 
\bigcap_{j=0}^i\mbox{\rm ker}({}_i\overline{\omega}c_{\hat{\theta_j}}) 
= \{0\}.    
$$ 
Hence, 
${}_i\overline{\omega} c_{(\theta_0, \ldots, \theta_i)}$
is injective.   
\par 
Since 
\begin{eqnarray*}
{} & { } & \dim_{\mathbb{K}}{\left(\frac{m_0^i\theta_0(2)}{m_0^{i+1}\theta_0(2)}\right)}
= 
2(i+1) \qquad \mbox{and}\\
{ } & { } & 
\dim_{\mathbb{K}} \left(
\frac{f^*m_0^i\theta_S(f)}{T\mathcal{R}_e(f)\cap f^*m_0^i\theta_S(f)+f^*m_0^{i+1}\theta_S(f)}
\right)
=2(i+1),
\end{eqnarray*}
where we have put $f=c_{(\theta_0, \ldots, \theta_i)}$,  
${}_i\overline{\omega}c_{(\theta_0, \ldots, \theta_i)}$ is an injective linear map between 
equidimensional vector spaces; this implies that it must be bijective.   
\hfill $\Box$ 
\par 
\medskip 
Theorem 2 is proved in Section 2.   
%%%%%%%%%%%%%%%%%%%%%%%%%%%%%%%%%%%%%%%%%%%%%%%%%%%%%%%%%%%%%%%%%%%%%%%% 
%%%%%%%%%%%%%%%%%%%%%%%%%%%%%%%%%%%%%%%%%%%%%%%%%%%%%%%%%%%%%%%%%%%%%%%% 
 %%%%%%%%%%%%%%%%%%%%%%%%%%%%%%%%%%%%%%%%%%%%%%%%%%%%%%%%%%%% 
\section{Proof of Theorem 2} 
Theorem 2 is proved by induction with respect to $i$.   
\par 
\medskip
\noindent 
\underline{\it The case $i=1$.}\quad     
Without loss of generality, we may assume that $\theta_0<\theta_1$.   
Then, it is sufficient to prove Theorem 2 for $R_{-\theta_0}\circ c_{(\theta_0,\theta_1)}$.    
Thus, in the following, we assume that $0=\theta_0<\theta_1<2\pi$ and $\theta_1\ne \pi$.   
\par 
It is easily seen that 
$$
\left\{
\left[\left(
\begin{array}{c}
2X \\ 
3Y
\end{array}
\right)
\right]
,\; 
\left[
\left(
\begin{array}{c}
2Y \\ 
0
\end{array}
\right)
\right]
\right\}
$$
is a basis of $\mbox{\rm ker}\left({}_1\overline{\omega}c_{0}\right)$. % (cf. \cite{nishimuraliftable}).
Therefore,  
{two vectors 
\begin{eqnarray*}
{ } & 
%\left\{
\left[
\left(
\begin{array}{c}
(2+\sin^2\theta_1)X-(\cos\theta_1\sin \theta_1)Y \\ 
-(\cos\theta_1\sin\theta_1)X+(2+\cos^2\theta_1)Y
\end{array}
\right)
\right]
,\; \\
{ } & 
\left[
\left(
\begin{array}{c}
-(2\cos\theta_1\sin\theta_1)X+(2\cos^2\theta_1)Y \\ 
-(2\sin^2\theta_1)X+(2\cos\theta_1\sin\theta_1)Y
\end{array}
\right)
\right]
%%\right\}
\end{eqnarray*}
constitute} 
a basis of $\mbox{\rm ker}\left({}_1\overline{\omega}c_{\theta_1}\right)$.   
\par 
Since 
$$
\dim_{\mathbb{K}}
\left(
\frac{m_0\theta_0(2)}{m_0^2\theta_0(2)}
\right)=4, 
$$
it is sufficient to show that 
$\mbox{\rm ker}\left({}_1\overline{\omega}c_0\right)\cap 
\mbox{\rm ker}\left({}_1\overline{\omega}c_{\theta_1}\right)
=\{0\}$.
In order to show this equality, it is sufficient to show that 
the determinant of the following matrix is not zero.   
$$
M=
\left(
{
\begin{array}{ccrc}
2 & 0 & 2+\sin^2\theta_1& -2\cos\theta_1\sin\theta_1 \\ 
0 & 2 & -\cos\theta_1\sin\theta_1 & 2\cos^2 \theta_1 \\
0 & 0 & -\cos\theta_1\sin\theta_1 & -2\sin^2\theta_1 \\ 
3 & 0 & 2+\cos^2\theta_1 & 2\cos\theta_1\sin\theta_1 
\end{array}
}
\right).
$$
Since $\det M={-20}\sin^2\theta_1$, it is not zero by the assumption of $\theta_1$.          
\par 
\medskip
\noindent 
%%%%%%%%%%%%%%%%%%%%%%%%%%%%%%%%%%%%%%%%%%%%%%%%%%%% 
\underline{\it The case $i=k+1$ under the assumption that Theorem 2 holds in the case $i=k$.}
\par 
\noindent 
%%%%%%%%%%%%%%%%%%%%%%%%%%%%%%%%%%%%%%%%%%%%%%%%%%%%     
 Since it has been assumed that Theorem 2 holds in the case $i=k$, Corollary 1 holds in the case $i=k$.      
{
It follows that for any  $j$ such that $0\le j\le k+1$, 
${}_{k+1}\overline{\omega}c_{\hat{\theta_j}}$ is surjective.    
%by the preparation the orem (for instance, see \cite{arnoldguseinzadevarchenko}).  
\begin{lemma}
For any $j$ such that $0\le j\le k+1$, the following holds:       
$$
\dim_{\mathbb{K}}\left(\mbox{\rm ker}\left({}_{k+1}\overline{\omega}c_{\hat{\theta_j}}\right)\right)=2.  
$$  
\end{lemma}
Although Lemma 1 is nothing but a special case of Proposition 4 in \cite{nishimuraliftable}, 
we give the proof of Lemma 1 here for the sake of clarity for the readers.      
\par 
\medskip
\noindent 
\underline{\it Proof of Lemma 1.}\qquad 
Take any $j$ $(0\le j\le k+1)$ and fix it.     
Put $f=c_{\hat{\theta_j}}$.       
For $f$, we need several notions defined in \cite{nishimura}.     
For any non-negative integer $\ell$, 
put ${}_\ell Q(f)=f^*m_0^\ell C_S/f^*m_0^{\ell+1}C_S$, 
${}_\ell \delta(f)=\dim_{\mathbb{K}}{}_\ell Q(f)$, and 
${}_\ell \gamma(f)=\dim_{\mathbb{K}}\mbox{ker}({}_\ell \overline{t}f)$, 
where ${}_\ell \overline{t}f: {}_\ell Q(f)\to {}_\ell Q(f)^2$ is defined by 
${}_\ell \overline{t}f([\eta])=[tf(\eta)]$.      
Then, it is easily seen that ${}_\ell \delta(f)=2(k +1)$ and 
${}_\ell \gamma(f)=1$.      
\par 
For $f$, $tf$ is injective and 
%${}_{k+1}\overline{\omega}c_{\hat{\theta_j}}$ is surjective.     
${}_{k+1}\overline{\omega}f$ is surjective.     
Therefore, 
\begin{eqnarray*}
{ } & { }  & 
\dim_{\mathbb{K}}\mbox{\rm ker}({}_{k+1}\overline{\omega}f) \\ 
{ } & = & 
\dim_{\mathbb{K}}\frac{m_0^{k+1}\theta_0(2)}{m_0^{k+2}\theta_0(2)} 
- 
\dim_{\mathbb{K}}\frac{f^*m_0^{k+1}\theta_S(f)}
{T\mathcal{R}_e(f)\cap f^*m_0^{k+1}\theta_S(f)+f^*m_0^{k+2}\theta_S(f)} \\ 
{ } & = & 
2(k+2) - 
\dim_{\mathbb{K}}\frac{\frac{f^*m_0^{k+1}\theta_S(f)}{f^*m_0^{k+2}\theta_S(f)}}
{\frac{T\mathcal{R}_e(f)\cap f^*m_0^{k+1}\theta_S(f)}
{T\mathcal{R}_e(f)\cap f^*m_0^{k+2}\theta_S(f)}} \\ 
{ } & = & 
2(k+2) - \left(2\cdot {}_{k+1}\delta(f)-\left({}_{k+1}\delta(f)-{}_{k+1}\gamma(f)+{}_k\gamma(f)\right)\right) \\ 
{ } & = & 
2(k+2) - 2(k+1) = 2.   
\end{eqnarray*} 
\hfill $\Box$ 
}
%Thus, by Theorem 1 and Proposition 1 of \cite{nishimuraliftable} (Proposition 1 of \cite{nishimuraliftable} 
%is an expressing formula of $\dim_{\mathbb{K}}{\rm ker}({}_{k+1}\overline{\omega}f)$ 
%by numerical invariants defined in \cite{nishimura} when ${}_{k+1}\overline{\omega}f$ is surjective), we have that 
%$$
%\dim_{\mathbb{K}}\left(\mbox{\rm ker}\left({}_{k+1}\overline{\omega}c_{\hat{\theta_j}}\right)\right)=2 
%\quad 
%(0\le j\le k+1). 
%$$   
\par 
{By Lemma 1, }
%On the other hand, we have that 
%$$
%\dim_{\mathbb{K}}\left(\frac{m_0^{k+1}\theta_0(2)}{m_0^{k+2}\theta_0(2)}\right)
%= 2(k+2).   
%$$
%Therefore, 
it is sufficient to show that 
$$
\frac{m_0^{k+1}\theta_0(2)}{m_0^{k+2}\theta_0(2)}=
\sum_{j=0}^{k+1}\mbox{\rm ker}\left({}_{k+1}\overline{\omega}c_{\hat{\theta_j}}\right).   
$$   
\par 
Let $(U, (X, Y))$ be a local coordinate system of $\mathbb{K}^2$ at the origin.   
Since $X^{k+1},$ $ X^kY, \ldots, Y^{k+1}$ are generators of $m_0^{k+1}$, 
for any $[\xi]\in \frac{m_0^{k+1}\theta_0(2)}{m_0^{k+2}\theta_0(2)}$, there must exist 
$\xi_1, \xi_2\in m_0^k\theta_0(2)$ such that 
$\xi=X\xi_1+Y\xi_2$.      
Then, by the assumption of induction we have the following:   
$$
[\xi_1], [\xi_2]\in 
\bigoplus_{j=0}^k\mbox{\rm ker}\left({}_k\overline{\omega}\left(c_{\hat{\theta_{k+1}}}\right)_{\hat{\theta_j}}\right).  
$$
Therefore, the following holds:   
$$
[\xi]\in 
\sum_{j=0}^k\mbox{\rm ker}\left({}_{k+1}\overline{\omega}\left(c_{\hat{\theta_{k+1}}}\right)_{\hat{\theta_j}}\right).  
$$
Hence, we may put 
$$
[\xi]=\sum_{j=0}^k[\widetilde{\xi}_j], 
$$
where $[\widetilde{\xi}_j]\in \mbox{\rm ker}\left({}_{k+1}\overline{\omega}\left(c_{\hat{\theta_{k+1}}}\right)_{\hat{\theta_j}}\right)$.   
Then, we have the following lemma.   
\begin{lemma}
For any $j$ such that $0\le j\le k$, the following holds:   
$$
[\widetilde{\xi}_j]\in 
\mbox{\rm ker}\left({}_{k+1}\overline{\omega}c_{\hat{\theta_j}}\right) +   
\mbox{\rm ker}\left({}_{k+1}\overline{\omega}c_{\hat{\theta_{k+1}}}\right). 
$$  
\end{lemma}
\par 
\medskip
\noindent 
\underline{\it Proof of Lemma 2.}\qquad 
Since $1\le k$, the set $\{0, 1, \ldots, k+1\}-\{j, k+1\}$ is not empty.    
Thus we may choose an element $\widetilde{\ell}\in \{0, 1, \ldots, k+1\}-\{j, k+1\}$.    
For $\widetilde{\ell}$, put $L=\{0, 1, \ldots, k+1\}-\{\widetilde{\ell}\}$.   
Then, from the assumption of induction, we have the following:   
$$
 \frac{m_0^k\theta_0(2)}{m_0^{k+1}\theta_0(2)}=
\bigoplus_{\ell\in L}\mbox{\rm ker}
\left({}_k\overline{\omega}\left(c_{\hat{\theta_{\widetilde{\ell}}}}\right)_{\hat{\theta_\ell}}\right).
$$
Hence, and since $[\widetilde{\xi}_j]\in \mbox{\rm ker}\left({}_{k+1}\overline{\omega}\left(c_{\hat{\theta_{k+1}}}\right)_{\hat{\theta_j}}\right)$, 
the following holds:   
$$
[\widetilde{\xi}_j]\in 
\left(
\sum_{\ell\in L}\mbox{\rm ker}\left({}_{k+1}\overline{\omega}
\left(c_{\hat{\theta_{\widetilde{\ell}}}}\right)_{\hat{\theta_\ell}}\right)
\right) 
\bigcap \mbox{\rm ker}\left({}_{k+1}\overline{\omega}\left(c_{\hat{\theta_{k+1}}}\right)_{\hat{\theta_j}}\right).
$$
Since $\widetilde{\ell}\ne j$ and $\widetilde{\ell}\ne k+1$, the intersection 
$$
\mbox{\rm ker}\left({}_{k+1}\overline{\omega}
\left(c_{\hat{\theta_{\widetilde{\ell}}}}\right)_{\hat{\theta_\ell}}\right)
\bigcap \mbox{\rm ker}\left({}_{k+1}\overline{\omega}\left(c_{\hat{\theta_{k+1}}}\right)_{\hat{\theta_j}}\right)
$$
is contained in 
$$
\mbox{\rm ker}\left({}_{k+1}\overline{\omega}c_{\hat{\theta_j}}\right) +   
\mbox{\rm ker}\left({}_{k+1}\overline{\omega}c_{\hat{\theta_{k+1}}}\right)
$$
for any $\ell\in L$.   
Thus, $[\widetilde{\xi}_j]$ must belong to  
$\mbox{\rm ker}\left({}_{k+1}\overline{\omega}c_{\hat{\theta_j}}\right)+ 
\mbox{\rm ker}\left({}_{k+1}\overline{\omega}c_{\hat{\theta_{k+1}}}\right)$.    
 \hfill $\Box$ 
\par 
\medskip 
By Lemma 2, we have the following:   
$$
[\xi]=\sum_{j=0}^k[\widetilde{\xi}_j]\in 
\sum_{m=0}^{k+1}\mbox{\rm ker}\left({}_{k+1}\overline{\omega}c_{\hat{\theta_m}}\right).
$$
Therefore, we have 
$$
\frac{m_0^{k+1}\theta_0(2)}{m_0^{k+2}\theta_0(2)}=
\sum_{m=0}^{k+1}\mbox{\rm ker}\left({}_{k+1}\overline{\omega}c_{\hat{\theta_m}}\right).   
$$ 
 \hfill $\Box$ 
 %%%%%%%%%%%%%%%%%%%%%%%%%%%%%%%%%%%%%%%%%%%%%%%%%%%%%%%%%%%%\]
%%%%%%%%%%%%%%%%%%%%%%%%%%%%%%%%%%%%%%%%%%%%%%%%%%%%%%%%%%%%%%%%%%%%%  

\bibliographystyle{amsplain}

\end{document}